\documentclass[12pt]{amsart}

\usepackage{amssymb}
\usepackage{mathrsfs}

\newtheorem{thm}{Theorem}[section]
\newtheorem{lem}[thm]{Lemma}

\theoremstyle{remark}
\newtheorem{remark}[thm]{Remark}

\newcommand{\symdif}{\triangle}
\newcommand\NULL{\emptyset}

\newcommand\xbf{\mathbf{x}}
\newcommand\ybf{\mathbf{y}}

\newcommand\comment[1]{}
\newcommand\zero{\mathtt{0}}
\newcommand\one{\mathtt{1}}
\newcommand\cat{{}^\smallfrown}
\newcommand\Bcal{\mathcal{B}}
\newcommand\Cbb{\mathbb{C}}
\newcommand\Nbb{\mathbb{N}}
\newcommand\Rbb{\mathbb{R}}
\newcommand\Tbb{U(1)}
\newcommand{\res}{\upharpoonright}

\title[Boolean action]{A Boolean action of $C(M,\Tbb)$ without\\ a spatial model and
a re-examination of the Cameron--Martin theorem}


 %should read ``2010'' but will not compile
\thanks{
We would like to thank Nathaniel Eldredge, Leonard Gross, and 
Alexander Kechris for reading an early draft of this article and offering 
a number of useful suggestions.
The research of the first author was supported by NSF grant DMS--0757507.
The research of the second author was supported by NSF grant DMS--1001623.
Any opinions, findings, and conclusions or recommendations expressed in this article
are those of the authors
and do not necessarily reflect the views of the NSF
}

\author[J. T. Moore]{Justin Tatch Moore}

\author[S. Solecki]{S{\l}awomir Solecki}

\address{Justin Tatch Moore:  Department of Mathematics \\ Cornell University\\
Ithaca, NY 14853--4201 \\ USA}

\email{{\tt justin@math.cornell.edu}}

\address{S{\l}awomir Solecki: Department of Mathematics \\
University of Illinois at Urbana-Champaign \\
1409 W. Green Street \\
Urbana, Illinois 61801-2975 \\ USA}

\email{{\tt ssolecki@math.uiuc.edu}}

\begin{document}

\begin{abstract}
We will demonstrate that if $M$ is an uncountable compact metric space,
then there is an action of the Polish group of all continuous functions
from $M$ to $\Tbb$ on a separable probability algebra which preserves
the measure and yet does not admit a point realization in the sense of Mackey.
This is achieved by exhibiting a strong form of ergodicity of the Boolean
action known as \emph{whirliness}.
This is in contrast with Mackey's point realization theorem,
which asserts that any measure preserving Boolean action of a
locally compact second countable group
on a separable probability algebra can be realized as
an action on the points of the associated probability space.
In the course of proving the main theorem, we will prove a
result concerning the infinite dimensional Gaussian measure space
$(\Rbb^\Nbb,\gamma_\infty)$ which is in contrast with the Cameron--Martin Theorem.
\end{abstract}

\maketitle

\section{Introduction}
Suppose that $G$ is a Polish group acting invariantly on a
probability space $(X,\mu)$ where $X$ is equipped with a Polish topology.\footnote{
Recall that \emph{Polish} is a synonym for separable and metrizable by a complete metric.}
Let $\Bcal_\mu$ denote the quotient of the Borel subsets of $X$ modulo the
$\mu$-null sets.
The original action induces a continuous action of $G$ on $\Bcal_\mu$, where
$\Bcal_\mu$ is topologized by the metric
$d_\mu (A,B) = \mu (A \symdif B)$.
An action of $G$ on $\Bcal_\mu$ induced by
a continuous homomorphism mapping $G$ into $\operatorname{Aut}(\Bcal_\mu,\mu)$
is commonly referred to as a \emph{Boolean action}.
In the situation above, the action of $G$ on
$(X,\mu)$ is commonly referred to as a \emph{spatial model}
(or \emph{point realization}) of this Boolean action.
In this paper, we will focus exclusively on Boolean actions on separable probability algebras;
the underlying set in a spatial model is by definition required to be a \emph{standard Borel space}
(i.e., the Borel sets arise from a Polish topology).

One can ask under what circumstances a Boolean action has a spatial model.
%Building on a result of von Neumann,
Mackey proved the first general result of this sort, demonstrating that if $G$ is a locally
compact Polish group, then every Boolean action of $G$ has a spatial model \cite{point_realize_trans}.
Mackey's proof relied heavily on the existence of Haar measure on $G$ and therefore
the situation for non-locally compact groups remained completely open.

Recently, Glasner, Tsirelson, and Weiss \cite{aut_grp_gauss_whirly}
introduced the notion of a \emph{whirly} Boolean action as a sufficient criterion for the 
\emph{non existence} of a spatial model for a Boolean action.
A Boolean action of a Polish group $G$ on $\Bcal_\mu$
is \emph{whirly} if whenever $V$ is a neighborhood of $1_G$ and $A,B \in \Bcal_\mu$
are positive, there is a $g$ in $V$ such that $\mu((g \cdot A) \cap B) > 0$.
It was
proved in \cite{aut_grp_gauss_whirly} that such Boolean actions not only do not have spatial models,
but in fact do not have non trivial factors with spatial models.
They also showed that all ergodic Boolean actions of Polish \emph{L{\'e}vy} groups are whirly, and exhibited such an
action for the L{\'e}vy group $L^0(\nu, \Tbb)$, where $\Tbb$ is the circle group and $\nu$ is a atomless Borel 
probability measure on a Polish space.
The construction of their Boolean action was generalized by
Pestov \cite{conc_whirl} to $L^0(\nu, S)$ for all non-trivial compact Polish groups $S$.
Here $L^0(\nu, S)$ is the group of all (equivalence classes) of $\nu$-measurable functions with values in $S$ taken
with the convergence in measure topology.

In the other direction, Glasner and Weiss proved that Boolean actions of
closed subgroups of the Polish group of all permutations of $\mathbb N$ always have spatial models
\cite{spatial_non-spatial}.
Kwiatkowska and the second author extended this and Mackey's results, proving that
Boolean actions of isometry groups of locally
compact separable metric spaces have spatial models \cite{spatial_isometries}.

It remained an open problem whether Mackey's theorem could be extended to groups of the form
$C(M,S)$ whenever $M$ is a compact metric space
and $S$ is a compact Polish group.
Here $C(M,S)$ is the Polish group of continuous functions from $M$ into $S$ with
the pointwise operation and the topology of uniform convergence. Part of the interest in such groups comes
from their position on the border between
the two classes of groups discussed above: the $L^0$ groups and the isometry groups of locally compact
separable metric spaces.
On the one hand, $C(M,S)$ maps continuously via inclusion onto a dense subgroup of
$L^0(\nu, S)$ for each Borel
probability measure $\nu$ on $M$.
On the other hand, the closed normal subgroups
$H\lhd C(M,S)$ for which $C(M,S)/H$ is compact separate the identity from non-identity elements
in $C(M,S)$.
This property places the groups $C(M,S)$ close to isometry groups of locally compact separable metric spaces,
since, by \cite[Theorem 1.2]{spatial_isometries}, its strengthening ``in each neighborhood of $1_G$ there is a closed normal subgroup $H$ for which $G/H$ is compact"
guarantees that a Polish group $G$ is the isometry group of a locally compact
separable metric space.

Kwiatkowska and the second author noticed, using the second author's spectral analysis of Boolean actions of $C(M, \Tbb)$,
that if $M$ is countable, then Boolean actions of $C(M,\Tbb)$ always admit spatial models.
In the present article, we will construct an example of an action of
$C(2^{\Nbb},\Tbb)$ on a separable atomless probability algebra which is whirly and hence
has no spatial model.
Notice that if $M$ is any uncountable compact metric space, then $M$ contains a homeomorphic
copy $P$ of the Cantor set.
Standard arguments show that the restriction map $f \mapsto f \restriction P$
defines a continuous epimorphism from $C(M,\Tbb)$ onto
$C(P,\Tbb) \simeq C(2^{\Nbb},\Tbb)$.
Furthermore, if $f:P \to \Tbb$,
then $f$ has a continuous extension to $M$ with the same distance from
the identity.
This yields the main result of the paper.

\begin{thm}\label{T:whirl}
If $M$ is an uncountable compact metric space, then
$C(M,\Tbb)$ admits a Boolean action on a separable atomless probability
algebra which is whirly.
\end{thm}
\noindent
The construction itself --- and Theorem \ref{main} which concerns it and is at the root
of the proof of Theorem~\ref{T:whirl} --- may
also be of independent interest from a probabilistic point of view.
Also, this appears to be the first whirly action whose whirliness does not
come from a form of concentration of measure on the acting group; see \cite{aut_grp_gauss_whirly},
\cite{spatial_non-spatial}, and \cite{conc_whirl}.

In the process of proving the main result, we will also establish
some results about the space $(\Rbb^\Nbb,\gamma_\infty)$, where $\gamma_\infty$
is the product measure arising from the standard Gaussian measure on $\Rbb$.

\begin{thm} \label{Xinfty_0-1_law}
If $a \in \Rbb$ and $K \subseteq \Rbb^{\Nbb}$ is a Borel set with $\gamma_\infty(K) > 0$,
then for $\gamma_\infty$-a.e. $y$ in $\Rbb^{\Nbb}$,
$\gamma_\infty ( \sqrt{1+a^2} K + a y ) > 0$.
\end{thm}

The result above is in contrast with
the Cameron-Martin Theorem (see, e.g., \cite[Example 2.3.5, Theorem 3.2.3]{gau_measure})
which implies that the translation of $\gamma_\infty$
by $x\in \Rbb^\Nbb$ is orthogonal to $\gamma_\infty$
unless $x$ is an element of $\ell^2$ (a set with $\gamma_\infty$-measure 0).
In fact, it is not difficult to show, using an analysis similar to that
in the proof of the Cameron-Martin Theorem, that if $(a,b) \ne (\pm1,0)$,
then the affine transformation $x \mapsto a x + b y$
sends a set of $\gamma_\infty$-measure $1$ to a set of
$\gamma_\infty$-measure $0$ for $\gamma_\infty$-a.e. $y$.
Also, Theorem \ref{Xinfty_0-1_law} does not hold for transformations
$x \mapsto a x + by$ unless $a^2 = b^2 + 1$.
This will be discussed in the final section.

Our notation is mostly standard.
We will take $\Nbb$ to include $0$ and all indexing will begin with $0$.
If $n$ is a natural number, then $2^n$ will be used to denote the collection of all binary
sequences of length $n$.
We will use $2^{<\omega}$ to denote the collection of all binary sequences of finite length.
If $\sigma$ is such a sequence and $n$ is not bigger than the length $|\sigma|$
of $\sigma$, then $\sigma\res n$ is the sequence
consisting of the first $n$ entries of $\sigma$.
If $\sigma$ and $\tau$ are sequences,
$\sigma \cat \tau$ will denote their concatenation.
If $x$ is a point in a metric space $X$ and $\delta > 0$, then
$N_\delta(x)$ denotes the open $\delta$-ball about $x$.
The variables $\epsilon$, $\delta$, and $s$ will always denote real numbers;
$i$, $j$, $k$, $m$, and $n$ will always denote integers.
We will use \cite{probability_measure} as our general
reference for probability.
Further information on measure theory and Polish group actions can be found in
\cite{geometry_quantum:Varadarajan}.

\section{The construction}

Let $\gamma$ denote the standard Gaussian measure on $\Cbb$ viewed as ${\mathbb R}^2$;
$\gamma$ has density
\[
\frac{1}{2\pi} e^{-\frac{x^2 + y^2}{2}}
\]
with respect to the Lebesgue measure on ${\mathbb R}^2$.
Define $I = 2^{<\omega} \cup \{*\}$ where $*$ is some symbol not in $2^{<\omega}$
and let $\gamma^I$ denote the product
measure on $\Cbb^I$.
Let $U_\sigma$ $(\sigma \in I)$ denote the coordinate functions on $\Cbb^I$,
regarded as complex random variables (i.e. these are i.i.d. standard complex Gaussian random variables).
Define random variables $Z_\sigma$ $(\sigma \in 2^{<\omega})$ recursively by
\[
Z_\NULL = U_*
\]
\[
Z_{\sigma \cat \zero} = \frac{Z_\sigma + U_\sigma}{\sqrt{2}}
\]
\[
Z_{\sigma \cat \one} = \frac{Z_\sigma - U_\sigma}{\sqrt{2}}.
\]
Observe that, for a given $\sigma$, $Z_{\sigma \cat \zero}$ and $Z_{\sigma \cat \one}$ are
independent.
It follows inductively that for each $n$,
$\{Z_\sigma : \sigma \in 2^n\} \cup \{U_\sigma : |\sigma| \geq n\}$
is a family of independent standard complex Gaussian random variables.
Hence if we define, for each $\sigma \in 2^n$ and $k \in \Nbb$,
\[
U_{\sigma,n+k} = \frac{1}{\sqrt{2}^k} \sum_{\tau \in 2^k} U_{\sigma \cat \tau},
\]
then 
$\{Z_\sigma : \sigma \in 2^n\} \cup \{U_{\sigma,k} : (\sigma \in 2^n) \land (k \geq n)\}$
is an independent family of standard complex Gaussian random variables. 

Now we define the set underlying the measure space on which the Boolean action 
of $C(2^\Nbb, U(1))$ will be performed. 
For each $n \in \Nbb$, define $X_n = \Cbb^{2^n}$ and let
$\gamma_n$ denote the product measure on $X_n$. 
If $n \in \Nbb$, define $\pi_{n+1,n} : X_{n+1} \to X_n$ by
\[
\pi_{n+1,n}(x)(\sigma) = \frac{x(\sigma \cat \zero) + x(\sigma \cat \one)}{\sqrt{2}}
\]
and if $n,k \in \Nbb$, let $\pi_{n+k,n} : X_{n+k} \to X_n$ be the composite projection
\[
\pi_{n+k,n} (x) (\sigma) =
\frac{1}{\sqrt{2}^k} \sum_{\tau \in 2^{k}} x (\sigma \cat \tau).
\]
It is easily verified that $\gamma_n(\pi_{n,m}^{-1}(A)) = \gamma_m(A)$ whenever
$A$ is a Borel subset of $X_m$ and $m < n$.
Define $X_\omega$ to be the collection of all functions $\xbf:2^{<\omega} \to \Cbb$ such that
for all $\sigma \in 2^{<\omega}$
\[
\xbf(\sigma) = \frac{\xbf(\sigma \cat \zero) + \xbf(\sigma \cat \one)}{\sqrt{2}}
\]
The projections $\pi_n:X_\omega \to X_n$ are defined by restriction.
Observe that $X_\omega$ is a closed linear subspace of $\Cbb^{2^{<\omega}}$ with
the product topology.

We define a measure $\gamma_\omega$ on $X_\omega$ as follows. Note that, for each $\sigma$,
\[
Z_\sigma = \frac{Z_{\sigma \cat \zero} + Z_{\sigma \cat \one}}{\sqrt{2}} 
\]
and, therefore, the random variables $Z_\sigma$ $(\sigma \in 2^{<\omega})$ induce a linear bijection from $\Cbb^I$ to $X_\omega$, 
where bijectivity follows from the identity 
\[
U_{\sigma} = \frac{Z_{\sigma \cat \zero} - Z_{\sigma \cat \one}}{\sqrt{2}}.
\]
Let $\Phi$ be the inverse of this bijection and
define $\gamma_\omega$ to be the pullback of $\gamma^I$ under $\Phi$:
$\gamma_\omega (A) = \gamma^I(\Phi(A))$.
It follows from the observations above that $\gamma_\omega$ satisfies
$\gamma_\omega(\pi_n^{-1}(A)) = \gamma_n(A)$
whenever $n \in \Nbb$ and $A \subseteq X_n$.
For each $\sigma \in 2^{<\omega}$ and $k\geq |\sigma|$, we will also regard
$Z_\sigma$, $U_\sigma$, and $U_{\sigma,k}$ as random variables in the measure
space $(X_\omega,\gamma_\omega)$ by identifying this space with $(\Cbb^I, \gamma^I)$ via 
$\Phi$.

Conditional probability will play a crucial role in our proof.
In particular, we will need the following result often known as the
Measure Disintegration Theorem.

\begin{thm}  (see \cite[5.14]{ergodic_theory:WE})
Let $(X,\mu)$ and $(Y,\nu)$ be standard probability spaces
and let $h:X \to Y$ be a Borel surjection such that
$\mu(h^{-1}(A)) = \nu(A)$ whenever $A$ is a Borel subset of $Y$.
Then there is a mapping $y \mapsto \mu_y$ defined on $Y$ which satisfies the following
conditions:
\begin{itemize}

\item for each $y$ in $Y$, $\mu_y$ is a Borel probability measure on $X$ and
$\mu_y (h^{-1}(y)) = 1$;

\item for each Borel $A \subseteq X$, $y \mapsto \mu_y (A)$ is a Borel
measurable function;

\item for every Borel $A \subseteq X$,
$\mu (A) = \int \mu_y (A) d \nu(y)$.
\end{itemize}
\end{thm}

For each $n$, fix a Borel measurable mapping
$x \mapsto \gamma_x$ defined on $X_n$ such that, for a given $x \in X^n$, $\gamma_x$ is a
Borel probability measure on $X_\omega$ with
\[
\gamma_x (\{\xbf \in X_\omega : \pi_n(\xbf) = x\}) = 1,
\]
\[
\gamma_\omega (A) = \int \gamma_x (A) d \gamma_n (x).
\]

Next we define the Boolean action of $C(2^{\Nbb},\Tbb)$ on $(X_\omega,\gamma_\omega)$.
Let $S_n$ denote all elements $g$ of $C(2^{\Nbb},\Tbb)$ such that
$g(x)$ depends only on the first $n$ coordinates of $x\in 2^\Nbb$.
Let $S_\omega$ denote the union of the sets $S_n$ as $n$ ranges over $\Nbb$.
Note that we can naturally identify $S_n$ with $\Tbb^{2^n}$.
If $m \geq n$, $g$ is in $S_n$, and $\sigma\in 2^m$,  we will use $g(\sigma)$
to denote the common value of $g(x)$ for those $x\in 2^\Nbb$ extending $\sigma$.
Observe that if $\xbf$ is a function defined on finite binary sequences of length
at least $n$ into $\Cbb$ and for every $\sigma$ of length at least $n$
\[
\xbf(\sigma) = \frac{\xbf(\sigma \cat \zero) + \xbf(\sigma \cat \one)}{\sqrt{2}},
\]
then there is a unique element of $X_\omega$ which extends $\xbf$.

If $\xbf$ is in $X_\omega$ and $g$ is in $S_n$,
then $g \cdot \xbf$ is defined by letting $(g\cdot \xbf)(\sigma)$, for $\sigma$ of
length at least $n$, be $g(\sigma)\xbf(\sigma)\in \Cbb$.
This is well defined for sequences $\sigma$ of length at least $n$ (notice that
$g$ is in $S_m$ for all $m \geq n$); the result extends uniquely to 
sequences of length less than $n$ to yield an element of $X_\omega$.
Thus, we have defined an action of $S_\omega$ on $X_\omega$ and hence on
$\Bcal_\omega$.
It should be clear that this action preserves $\gamma_\omega$.

\begin{lem} The action of $S_\omega$ on ${\mathcal B}_\omega$ extends to a continuous 
action of $C(2^{\mathbb N}, U(1))$ on ${\mathcal B}_\omega$. 
\end{lem}

\begin{proof} The lemma is a consequence of the following statement. 

For every $W \subseteq X_k$ which is a product of open disks
and for every $\epsilon > 0$,
there is a $\delta > 0$ such that if $n >k$ and
$g,h \in S_n$ with $|g(\sigma)-h(\sigma)| < \delta$ for
all $\sigma \in 2^n$, then
$\gamma_n ((g \cdot \pi_{n,k}^{-1}(W)) \symdif (h \cdot \pi_{n,k}^{-1}(W))) < \epsilon$. 

To prove the above statement, first observe that it suffices to prove it for $k=0$ 
(that is, $W$ is an open disk in $\Cbb$).
Let $z$ and $r$ be such that $W = N_r(z)$.
Let $s >0$ be sufficiently small that
\[
\gamma (N_{r+s}(x) \setminus N_{r-s} (x)) < \epsilon/3.
\]
Define $\delta = s \sqrt{\epsilon/3}$.
Now suppose that $g,h \in S_n$ are given with $|g(\sigma)-h(\sigma)|<\delta$ for all $\sigma\in 2^n$.
Notice that
\[
x \mapsto \pi_{n,0}(g^{-1} \cdot x) - \pi_{n,0}(h^{-1} \cdot x)\in \Cbb
\]
defines a complex Gaussian random variable on $(X_n,\gamma_n)$
with mean $0$ and variance less than $\delta^2$.
In particular, by Chebyshev's inequality
\[
\gamma_n (\{x \in X_n : |\pi_{n,0}(g^{-1} \cdot x) - \pi_{n,0}(h^{-1} \cdot x)| \geq s\}) \leq \frac{\epsilon}{3}.
\]
Finally, observe that $x \in g \cdot \pi_{n,0}^{-1}(W)$ is equivalent to
$\pi_{n,0}(g^{-1} \cdot x) \in W$ (and similarly for $h$).
Thus, if we let 
\[
E = (g \cdot \pi_{n,0}^{-1}(W)) \symdif (h \cdot \pi_{n,0}^{-1}(W))
\]
\[
E_0 = \{x \in X_n : \pi_{n,0} (x) \in N_{r+s}(x) \setminus N_{r-s} (x)\}
\]
\[
E_1 = \{x \in X_n : |\pi_{n,0}(g^{-1} \cdot x) - \pi_{n,0}(h^{-1} \cdot x)| \geq s\}, 
\]
then $E \subseteq gE_0 \cup hE_0 \cup E_1$ and hence $\gamma_n(E) < \epsilon$. The statement 
and the lemma follow.
\end{proof}

We finish this section by recalling some
standard facts about the space $(\Cbb^I,\gamma^I)$.
Let $H$ denote the collection of all elements of $\Cbb^I$ with finite support, regarded
as a group under addition.

\begin{lem} \label{quasi-inv}
For every Borel $A \subseteq \Cbb^I$ and every $h \in H$,
$\gamma^I(A) > 0$ if and only if $\gamma^I(A + h) > 0$.
\end{lem}

\begin{proof}
This follows from the fact that translations are
absolutely continuous in finite dimensional powers of $(\Cbb,\gamma)$
and the properties of product measures.
(It is also a trivial consequence of the Cameron-Martin Theorem; see \cite[Example 2.3.5, Theorem 3.2.3]{gau_measure}.)
\end{proof}

\begin{lem} \label{0-1_law}
If $A \subseteq \Cbb^I$ is Borel and $H$-invariant (i.e., $A + H = A$),
then $\gamma^I(A)$ is either $0$ or $1$.
\end{lem}

\begin{proof}
This is a special case of Kolmogorov's Zero-One Law.
\end{proof}

\section{Proof of the main theorem}

For each $s \in \Rbb$, $k \in \Nbb$ define the following element of $C(2^{\Nbb},\Tbb)$:
\[
g_{s,k}(x) =
\frac{1 + i s (-1)^{x(k)}}{\sqrt{1 + s^2}}
\]

\begin{lem} \label{indep_rand}
For each $n$ and $\gamma_n$-a.e. $z\in X_n$,
\[
\{U_{\sigma,k} : (\sigma \in 2^{n}) \land (n\leq k)\}
\]
is a family of mutually independent standard complex Gaussian random variables with respect to $\gamma_z$.
\end{lem}

\begin{proof}
This follows immediately from the observation made above that
\[
\{Z_\sigma : \sigma \in 2^n\} \cup \{U_{\sigma,k} : (\sigma \in 2^n) \land (n \leq k)\}
\]
is an independent family of standard complex Gaussian random variables with respect to $\gamma_\omega$.
\end{proof}

\begin{lem} \label{cond_indep_lem}
If $n \in \Nbb$ and $K \subseteq X_n$ is Borel with $\gamma_n(K) > 0$,
then for each $s > 0$ and $\gamma_n$-a.e. $z \in X_n$
\[
g_{s,k}\cdot \pi_n^{-1}(K),\;\hbox{ for } n \leq k
\]
are mutually independent with respect to $\gamma_z$,
each having $\gamma_z$-measure
$\gamma_n \Big( \frac{\sqrt{1 + s^2} K - z}{s} \Big)$.
\end{lem}

\begin{proof}
Let $U_{n,k}$ denote the vector valued random variable
\[
(U_{\sigma,k} : \sigma \in 2^n ).
\]
Let $\ybf$ be an element of $X_\omega$ and $\sigma$ be an element of $2^n$.
Observe that
the value of $\pi_n(g_{s,n} \cdot \ybf)$ at $\sigma$ is
\[
\frac{\ybf(\sigma \cat \zero) + i s  \ybf(\sigma \cat \zero) + \ybf(\sigma \cat \one) - i s \ybf(\sigma \cat \one)}
{\sqrt{2}\sqrt{1 + s^2}}
\]
\[
= \frac{1}{\sqrt{1 + s^2}}
\Big(
\frac{\ybf(\sigma \cat \zero) + \ybf(\sigma \cat \one)}{\sqrt{2}} + 
i s \frac{ \ybf(\sigma \cat \zero) - \ybf(\sigma \cat \one)}{\sqrt{2}}
\Big)
\]
and thus
\[
\pi_n(g_{s,n} \cdot \ybf) = \frac{\pi_n(\ybf) + i s U_{n,n}(\ybf)}{\sqrt{1 + s^2}}
\]
More generally, if $n \leq k$, then
\[
\pi_n(g_{s,k} \cdot \ybf)  = \frac{\pi_n(\ybf) + i s U_{n,k}(\ybf)}{\sqrt{1 + s^2}}.
\]

Fix $z\in X_n$. It follows that if $\pi_n(\ybf)=z$, then
$\pi_n(g_{s,k} \cdot \ybf)$ is in $K$ if and only if
\[
\frac{1}{\sqrt{1 + s^2}} (z + i s U_{n,k}(\ybf)) \in K
\]
if and only if
\[
U_{n,k}(\ybf) \in -i\frac{ \sqrt{1 + s^2} K - z}{s}.
\]
The conclusion follows from Lemma \ref{indep_rand} after
observing that multiplication by $-i$ preserves the standard Gaussian measure.
\end{proof}

\begin{lem} \label{convolve}
If $K \subseteq X_\omega$ is Borel and $a \in \Rbb$, then
\[
\int \gamma_\omega ( \sqrt{1 + a^2} K + a \ybf ) d \gamma_\omega (\ybf) = \gamma_\omega(K)
\]
\end{lem}

\begin{remark} \label{convolve_rem}
Note that, since there is a continuous linear measure preserving bijection between
$(X_\omega,\gamma_\omega)$ and $(\Rbb^\Nbb,\gamma_\infty)$, an analogous identity holds
for the measure space $(\Rbb^\Nbb,\gamma_\infty)$.
This can also be verified by direct computation.
\end{remark}

\begin{proof} We first show that 
if $n \in \Nbb$, $K \subseteq X_n$ is Borel, and $a \in \Rbb$,
then
\begin{equation}\notag
\int \gamma_n (\sqrt{1 + a^2} K + a z ) d \gamma_n (z) = \gamma_n(K)
\end{equation}
If $a =0$, the above equality clearly holds. If $a \ne 0$, set $s = -1/a$.
Applying Lemma \ref{cond_indep_lem} we have
\[
\gamma_n (K) = \gamma_\omega(\pi_n^{-1} (K)) =
\gamma_\omega (g_{s,n} \cdot \pi_n^{-1} (K))
\]
\[
= \int \gamma_z (g_{s,n} \cdot \pi_n^{-1} (K)) d \gamma_n (z)
\]
\[
= \int \gamma_n ( \sqrt{1 + a^2} K + az ) d \gamma_n(z),
\]
as required. 

Observe now that
$\nu (K) = \int \gamma_\omega ( \sqrt{1 + a^2} K + a \ybf ) d \gamma_\omega (\ybf)$
defines a Borel probability measure on $X_\omega$.
By the equality proved above, this measure coincides with $\gamma_\omega$ on
sets of the form $\pi_n^{-1}(K)$ for $K$ a Borel subset of $X_n$.
Therefore $\nu$ equals $\gamma_\omega$.
% \subseteq 
%Since $\gamma_\omega$ is a Radon measure, it suffices to prove the lemma
%for compact $K$.
%Applying Fact \ref{cpt_approx}, Fact \ref{pi_linear}, the Monotone Convergence Theorem,
%Fact \ref{restricted_dependence},  Lemma \ref{fin_dim_convolv}, and Fact \ref{cpt_approx}
%in sequence, we obtain:
%\[
%\int \gamma_\omega ( \sqrt{1 + a^2} K + a \ybf ) d \gamma_\omega(\ybf) =
%\int \lim_{n \to \infty} \gamma_n ( \pi_n (\sqrt{1 + a^2}K + a\ybf )) d \gamma_\omega (\ybf)
%\]
%\[
%= \int \lim_{n \to \infty} \gamma_n (\sqrt{1 + a^2}\pi_n(K) + a \pi_n(\ybf) ) d \gamma_\omega(\ybf)
%\]
%\[
%= \lim_{n \to \infty} \int \gamma_n (\sqrt{1 + a^2}\pi_n(K) + a \pi_n(\ybf) ) d \gamma_\omega (\ybf)
%\]
%\[
%= \lim_{n \to \infty} \int \gamma_n (\sqrt{1 + a^2}\pi_n(K) + a y ) d \gamma_n(y)
%\]
%\[
%= \lim_{n \to \infty} \gamma_n(\pi_n(K)) = \gamma_\omega(K).
%\]
\end{proof}

\begin{lem} \label{0-1_lem}
If $a \in \Rbb$ and $K \subseteq X_\omega$ is Borel with $\gamma_\omega(K) > 0$, then
\[
\gamma_\omega(\{ \xbf \in X_\omega :
\gamma_\omega (\sqrt{1 + a^2} K + a \xbf ) > 0
\}) = 1
\]
\end{lem}

\begin{remark}
Theorem \ref{Xinfty_0-1_law} follows immediately by again
making the observation that there is a
continuous linear measure preserving bijection between
$(X_\omega,\gamma_\omega)$ and $(\Rbb^\Nbb,\gamma_\infty)$.
\end{remark}

\begin{proof}
Let $K$ and $a$ be given as in the statement of the lemma and define
\[
A = \{ \xbf \in X_\omega :
\gamma_\omega (\sqrt{1 + a^2} K + a \xbf ) > 0
\}
\]
and let $B$ be the image of $A$ under $\Phi$.
First observe that by Lemma \ref{convolve}, $\gamma_\omega(A) > 0$.
If $\xbf$ is in $A$, then the linearity of $\Phi$ implies that
\[
\gamma^I ( \sqrt{1+a^2}\Phi(K) + a \Phi(\xbf) ) > 0.
\]
By Lemma \ref{quasi-inv}, if $h$ is in $H$, then
\[
\gamma^I ( \sqrt{1+a^2}\Phi(K) + a \Phi(\xbf) + a h ) > 0
\]
and therefore $\xbf + \Phi^{-1}(h)$ is also in $A$.
This means that $B$ is $H$-invariant and therefore must have measure $1$ by Lemma~\ref{0-1_law}.
Since $\gamma_\omega(A) = \gamma^I(B)$, this completes the proof.
\end{proof}

We are now ready to establish that our near action is whirly.

\begin{thm} \label{main}
For every Borel $K \subseteq X_\omega$ with $0 < \gamma_\omega(K)$ and
for every $\epsilon > 0$ there are $m$ and $n$ with
\[
\gamma_\omega \Big( \bigcup_{k=n}^{n+m} g_{\epsilon,k} \cdot K \Big) > 1 - \epsilon.
\]
\end{thm}

\begin{remark}
The functions $g_{\epsilon,k}$ are in fact Lipschitz.
Notice, however, that while the uniform distance from $g_{\epsilon,k}$ to $1$ does
not depend on $k$ (and tends to $0$ as $\epsilon$ approaches $0$),
the Lipschitz norm of $g_{\epsilon,k}$ tends to infinity with respect
to $k$ for any fixed $\epsilon >0$.
\end{remark}

\begin{proof}
Let $K$ and $\epsilon$ be given as in the statement of the theorem. We can, and do, assume that $K$ is compact. 
Set $a = -1/\epsilon$ and fix a $\delta > 0$ such that
\[
\gamma_\omega (\{\ybf \in X_\omega : \gamma_\omega (\sqrt{1 + a^2} K + a \ybf ) \geq \delta\})
> \sqrt{1- \epsilon/2}.
\]
This is possible by Lemma \ref{0-1_lem}.
Let $m$ be sufficiently large such that $1- (1-\delta)^m \geq \sqrt{1-\epsilon/2}$.
Next by applying compactness of $K$, fix an $n$ which is sufficiently large such that
\[
\gamma_\omega(K) > \gamma_n(\pi_n(K)) - \frac{\epsilon}{2m}.
\]

Define $K' = \pi_n(K)$ and set
\[
A = \{z \in X_n : \gamma_n (\sqrt{1 + a^2} K' + a z ) \geq \delta\}
\]
noting that $\gamma_n(A) > \sqrt{1-\epsilon/2}$.
By Lemma \ref{cond_indep_lem}, for each $z$ in $A$, we have that the sets
$g_{\epsilon,k} \cdot \pi_n^{-1}(K')$, $n \leq k < n + m$, are $\gamma_z$-mutually independent with $\gamma_z$-measure
at least $\delta$.
Their union, therefore, has $\gamma_z$-measure at least $1-(1-\delta)^m \geq \sqrt{1-\epsilon/2}$.
Integrating, we obtain
\[
\gamma_\omega \Big(\bigcup_{k=n}^{n+m-1} g_{\epsilon,k} \cdot \pi_n^{-1}(K') \Big) =
\int_{X_n} \gamma_z \Big(\bigcup_{k=n}^{n+m-1} g_{\epsilon,k} \cdot \pi_n^{-1}(K') \Big) d \gamma_n(z)
\]
\[
\geq \int_A \gamma_z \Big(\bigcup_{k=n}^{n+m-1} g_{\epsilon,k} \cdot \pi_n^{-1}(K') \Big) d \gamma_n(z)
> \Big(\sqrt{1-\frac{\epsilon}{2}}\Big)\Big(\sqrt{1-\frac{\epsilon}{2}}\Big) = 1-\frac{\epsilon}{2}
\]
Finally,
\[
\gamma_\omega \Big(\bigcup_{k=n}^{n+m-1} g_{\epsilon,k} \cdot K \Big) >
\gamma_\omega \Big(\bigcup_{k=n}^{n+m-1} g_{\epsilon,k} \cdot \pi_n^{-1}(K') \Big) - m \frac{\epsilon}{2m}
>  1-\epsilon.
\]
\end{proof}

\section{Concluding remarks}

It is natural to ask whether the main construction of this paper is possible if we replace
$C(M,\Tbb)$ by $C^\infty(M,\Tbb)$, where $M$ is now some smooth manifold and
$C^\infty(M,\Tbb)$ is equipped with the metric arising from the $C^\infty$-norm.
Suppose that $M$ is given and $P \subseteq M$ is homeomorphic to $2^\Nbb$.
The Boolean action itself can be constructed exactly as above.
The key difference, however, is that for a fixed $\epsilon > 0$,
the extensions of the functions $g_{\epsilon,k}$ to $M$ will
necessarily have norms which diverge to infinity as $k \to \infty$.
The role of the functions $g_{\epsilon,k}$
seems so essential to our argument that it is tempting to conjecture that 
if $M$ is a smooth manifold and $S$ is a compact Lie group, then
every Boolean action of $C^\infty (M,S)$ on a separable probability algebra
has a spatial model.

We will finish by mentioning that
Theorem \ref{Xinfty_0-1_law} is sharp in the sense that
it only applies to the affine transformations
$x \mapsto a x + by$ when $a^2 = b^2 + 1$.
To see this define, for $x \in \Rbb^\Nbb$,
\[
s(x) = \lim_{n \to \infty} \left(\frac{1}{n} \sum_{k < n} |x (k)|^2 \right)^{\frac{1}{2}}.
\]
While this limit may not exists,
$s (x) = 1$ for $\gamma_\infty$-a.e. $x$; set
\[
K = \{x \in \Rbb^\Nbb : s(x) = 1\}.
\]
Furthermore, if $x$ and $y$ are in $\Rbb^{\Nbb}$ and $a$ is in $\Rbb$,
then $s(a x) = a s(x)$.
Next observe that if $a,b \in \Rbb$,
then by Remark \ref{convolve_rem},
\[
\gamma_\infty ( \frac{a}{\sqrt{1+ b^2}} K) = \int \gamma_\infty (\sqrt{1 + b^2} \frac{a}{\sqrt{1 + b^2}} K + b y) d \gamma_\infty (y)
\]
\[
=\int \gamma_\infty (a K + b y) d \gamma_\infty (y).
\]
If $a^2 \ne b^2 + 1$, then
$(\frac{a}{\sqrt{1 + b^2}} K) \cap K = \emptyset$.
In particular, $\gamma_\infty (\frac{a}{\sqrt{1 + b^2}} K) = 0$ and hence
$\gamma_\infty (a K + b y) = 0$ for $\gamma_\infty$-a.e. $y$.

\def\Dbar{\leavevmode\lower.6ex\hbox to 0pt{\hskip-.23ex \accent"16\hss}D}

\end{document}